%Paper: alg-geom/9312002
%From: James Madden <madden@marais.math.lsu.edu>
%Date: Fri, 3 Dec 93 09:24:56 -0600

%TeX version 3.1
\magnification=\magstep1
\mathsurround=1pt

\outer\def\myproclaim #1: #2\par{\medbreak\vskip-\parskip
 {\bf#1:\enspace}{\sl#2}\par
  \ifdim\lastskip<\medskipamount

\removelastskip\penalty55\medskip\fi}

\def\a{\alpha}
\def\b{\beta}
\def\g{\gamma}

\def\cntr{\mathop{\rm cntr}}
\def\dim{\mathop{\rm dim}}

\def\m{{\bf m}}
\def\ord{\mathop{\rm ord}}

\def\p{{\bf p}}
\def\q{{\bf q}}

\def\res{\mathord{\rm res}}
\def\sepab{\langle \alpha , \beta \rangle}
\def\sepag{\langle \alpha , \gamma \rangle}

\def\sper{\mathop {\rm Spec}_r }

\def\supp{\mathop{\rm supp}}
\def\td{\mathop{\rm tr.deg.}}
\def\qed{\hfil///// \parfillskip0pt\par\smallskip\parfillskip0pt
plus1fil}
\def\R{{\bf R}}
\def\Q{{\bf Q}}
\def\Z{{\bf Z}}

\centerline{\bf  Complete ideals defined by sign conditions}
\centerline{\bf  and the real spectrum of a two-dimensional local
ring}
\vskip 1pc
\centerline{Dean Alvis}
\centerline{Bernard L. Johnston}
\centerline{James J. Madden}

\bigskip

{\bf 0. Introduction.}  Let $(A,\m)$ be a regular local ring and let
$\a$ and $\b$ be points of the real spectrum of $A$ centered at $\m$.
$\a$ and $\b$ may be viewed as total orderings on quotients of $A$.
Associated with $\a$ and $\b$, there is the so-called ``separating
ideal'', $\sepab\subseteq A$, which is generated by all $a\in A$ such
that $a$ is non-negative with respect to $\a$ and $-a$ is
non-negative with respect to $\b$.  $\sepab$ is a valuation ideal for
the valuation canonically associated with $\a$ (or $\b$), and hence
is complete.  It is known that a thorough understanding of $\sepab$
would contribute greatly to a solution of the long-standing
Pierce-Birkhoff conjecture (see [M]) but up to now no good techniques
for working with it have been known.  In this paper, we investigate
$\sepab$ by applying quadratic transforms to $(A,\m)$ and using
Zariski's theory of complete ideals in two-dimensional regular local
rings (see [Z]) to analyze how $\sepab$ is affected.    Suppose that
$(A', \m')$ is a quadratic transform of $A$.  Under natural
hypotheses, $\a$ and $\b$ induce points $\a'$ and $\b'$ in the real
spectrum of $A'$ which are centered at $\m'$.  Our main result,
Theorem 4.7, is a formula which relates $\langle \a', \b' \rangle
\subseteq A'$ with the ideal transform $\sepab'$ of $\sepab$ in $A'$.
It says that if $A$ is two-dimensional and has real closed residue
field, and if $\sepab$ is not the maximal ideal, then $\langle \a',
\b' \rangle = \sepab'$.  Applications of this result are presented in
detail elsewhere (see [MR] and [MS]).  The applications show that the
transformation formula provides an essentially complete understanding
of separating ideals in two dimensional regular algebras over real
closed fields.

Here is a summary of the contents.  In Section 1,  we review notation
for valuations and for the real spectrum, and then make some
observations on valuations induced by points of the real spectrum.
In Section 2, we  prove (Proposition 2.2) that separating ideals are
simple if a certain technical condition is satisfied.  In section 3,
we examine the behavior of the real spectrum under quadratic
transformation, and we prove some general facts about the effect of
quadratic transformations on separating ideals.  We also give an
example which indicates that in dimensions 3 and higher the behavior
of separating ideals under quadratic transformations will be very
difficult to analyze.  In section 4, we consider 2-dimensional
regular local rings, and we prove the transformation formula
mentioned above.  With Zariski's theory and the results of sections 1
through 3 at hand, the hardest part of the proof of the
transformation formula is Theorem 4.4.  A sampling of some of the
results which will be presented in [MR] and [MS] is given in the last
section.

The present paper had a rather lengthly gestation.  Alvis and Madden
worked together on the type of ideals studied here in 1989, but were
unaware of Zariski's work on complete ideals at that time.  Alvis
wrote several computer programs which searched for generators for
separating ideals.  Without the wealth of examples found this way,
connections to quadratic transforms would not have been noticed later
on.  Madden learned about complete ideals in two-dimensional rings
from Johnston in 1990-91, and together they conjectured a version of
Theorem 4.7 in early 1991.  The proof of 4.7 was completed by Madden
in summers of 1991 and 1992, with the support of NSF-DMS 9104427.

\bigskip
{\bf 1. Valuations and orderings.}  We describe first the notation we
use when discussing valuations.  Let  $A$ be a noetherian ring and
let $v$ be an additive valuation.  In other words, $v$ is a map from
$A$ to $\Gamma \cup \{+\infty\}$, where $\Gamma$  is a totally
ordered abelian group, and $v$ satisfies  $v(a+b) \geq \min\{\,v(a),
v(b)\,\}$, $v(ab) = v(a)+v(b)$, $v(0)=+\infty$ and $v(1)=0$.  An
ideal $I \subseteq A$ is a {\it v-ideal} if  $a \in I$ and ${\it
v}(b) \geq {\it v}(a)$ imply $b\in I$.  The collection of all {\it
v}-ideals in $A$ is totally ordered by inclusion.  The smallest {\it
v}-ideal in $A$  (namely, $v^{-1}(+\infty)$) is called the {\it
support} of  $v$ in $A$ and is denoted $\supp_A v$, or just $\supp v$
if no confusion is possible.  The largest proper {\it v}-ideal in $A$
is called the {\it center} of  $v$ in $A$, and is denoted $\cntr_A
v$, or just $\cntr v$.  We say that $v$ is {\it non-trivial} if
$\cntr v \not= \supp v$.  The largest {\it v}-ideal properly
contained in the {\it v}-ideal $I$  ---this exists because $A$ is
noetherian---is called the {\it successor of} $I$ and is denoted
$I^{\hbox{\it v}}$.  The valuation $v$ induces a valuation $\overline
{v}$ on the fraction field of $A/\supp v$.  The valuation ring of
this valuation is denoted $V_v$ and its maximal ideal is denoted
$\m_v$.  The residue field of $v$ is $K_v : = V_v/\m_v$.  Let
$d_v:=A/\cntr v$, and let $k_v$ denote th fraction field of $d_v$.
We have the containments $d_v\subseteq k_v \subseteq K_v$.  The {\it
dimension} of $v$ is the transcendence degree of $K_v$ over $k_v$,
$\dim v:=\td (K_v/k_v)$.

We now describe the notation we use when dealing with the real
spectrum.  It is close to [BCR], Chapter 7.  $\sper A$ is the set of
prime cones in $A$.  Each $\a \in \sper A$ determines a pair
$(\supp\a, \leq_\a)$, where $\supp\a = \a \cap -\a$ is a real prime
ideal of $A$ and $\leq_\a$ is a total ordering of $A/\supp\a$
satisfying the usual conditions of compatibility with the ring
structure.  Given a real prime ideal $I$ and a total order $\leq_0$
on $A/I$, there is a prime cone $\a = \{\, a\in A \mid 0 \leq_0 a+I
\,\}$, with $\supp\a = I$ and $\leq_\a = \leq_0$.  We let $A[\a ]$
denote $A/\supp\a$ endowed with the ordering $\leq_\a$.  If $a \in
A$, then $a(\a)$ denotes the residue of $a$ in $A[\a ]$.  Also,
$A(\a)$ denotes the ordered fraction field of $A[\a ]$.  If
$\a\supset \b$, we say $\a$ is a specialization of $\b$.  In this
case, there is an order-preserving homomorphism $A[\b ]\to A[\a ]$.

There is a valuation on  $A$ naturally associated with $\a$. Let
$$V_\a:=\{\, a \in A(\a) \mid \exists x \in A[\a ] \enskip |a| \leq
|x| \,\}.$$   This is a valuation ring in $A(\a)$.  Let
$\overline{v_\a}:A(\a)\rightarrow\Gamma_\a \cup \{+\infty\}$ be the
valuation associated with $V_\a$, and let $v_\a :A \to \Gamma_\a \cup
\{+\infty\}$ be the valuation obtained by mapping $A$ to $A(\a)$ and
then applying $\overline {v_\a}$.  Note that $v_\a$ is non-negative
on $A$.  It seems worth explicit mention that $v_\a$ may be trivial
even in cases where the ordering $\a$ is in no sense trivial, {\it
viz}., when $A$ is a field.  This is also possible when $A$ is not a
field, for example, let $A=\Q[x_1, \dots , x_n]$ and let $\a$ be
induced by the embedding of $A$ in $\R$ which identifies the
indeterminates with $n$ algebraically independent real transcendental
numbers.

The $v_\a$-ideals have a nice description in terms of the ordering of
$A[\a]$.  An ideal $I$ in $A[\a ]$ is {\it convex} if $a\in I$ and
$|b|\leq |a|$ implies $b\in I$.  An ideal of $A$ is called an
$\a$-ideal if it is the pull-back of a convex ideal of $A[\a ]$.

\smallskip
\myproclaim Lemma 1.1: $I \subseteq A$ is an $\a$-ideal if and only
if it is a $v_\a$-ideal.

\smallskip
{\it Proof}.  For all

$a, b\in A[\a ]$, $0 \not= a$, we have the following equivalences:
$$\displaylines{
\quad\exists x\in A[\a ]\enskip |b|\leq |xa| \Leftrightarrow \exists
x\in A[\a ]\enskip|b|/|a| \leq |x|\hfill\cr
\hfill{}\Leftrightarrow |b|/|a|\in V_\a \Leftrightarrow
\overline{v_\a} (|b|/|a|)\geq 0 \Leftrightarrow \overline{v_\a} (b)
\geq \overline{v_\a} (a).\quad\cr
}$$
Suppose

$I \subseteq A[\a]$ is a $\overline{v_\a}$-ideal.  If $a\in I$, $0
\not= a$, and $|b|\leq |a|$, then $\overline{v_\a} (b) \geq
\overline{v_\a} (a)$, so $b\in I$.  Thus, $I$ is convex.

On the other hand, if $I$ is convex and $a\in I$ and

$\overline{v_\a} (b) \geq \overline{v_\a} (a)$, then $|b|$ is less
than or equal to some element of $I$, so $I$ is a
$\overline{v_\a}$-ideal.
\qed

Some important concepts related to $v_\a$ can be defined in terms of
the ordering $\a$.  For example, $\cntr_Av_\a$ is the largest proper
$\a$-ideal of $A$.  We shall denote this by $\cntr\a$.  Let $\m_\a$
be the maximal ideal of $V_\a$ and let $K_\a: = V_\a/\m_\a$.  Note
that $K_\a$ is naturally an ordered field.   Let $f$ be the map $A
\to V_\a \to K_\a$.  The kernel of $f$ is $\cntr\a$.    The {\sl
image} of $f$ is a sub-domain  $d_\a \subseteq  K_\a$.  This notation
is summarized by the following diagram, in which the rows are exact
the vertical arrows are inclusions:
$$\matrix{0 &\to &\cntr\a/\supp\a&\to &A[\a]    &\to &d_\a &\to &0\cr
	    &    &\downarrow& &\downarrow&&\downarrow&&\cr
	  0 &\to &\m_\a  &\to &V_\a &\to &K_\a &\to
                                                  &0\cr}.$$
The fraction field of $d_\a$ is denoted $k_\a$; $\dim\a $ is defined
to be  $\td(K_\a/k_\a)$.  (Caution: $\dim\a$ should not be confused
with $\dim\supp\a$, by which is meant the Krull dimension of
$A/\supp\a$.)

Since every element of $V_\a$ is between two elements of  $A[\a]$,
every element of $K_\a$ is between two elements of $d_\a$, {\it
i.e.}, $d_\a$ is cofinal in $K_\a$.  This has an important
consequence when $d_\a$ happens to be the field of real numbers (as,
for example, when $A$ is an $\R$-algebra).  In this case, $d_\a =
K_\a$, since $\R$ is not cofinal in any ordered field which is a
proper extension.  In general, of course, one must deal with
valuations centered at maximal ideals which are not zero-dimensional,
{\it e.g.}, the order valuation on a local ring.

\smallskip
{\bf Examples 1.}  We first describe a convenient method of producing
examples, which will be used repeatedly.  Suppose $B$ is a
totally-ordered domain and $\phi: A \to B$ is any ring homomorphism.
Let $\leq_\phi$ be the ordering on $A/{\rm ker}\phi$ induced by its
inclusion in $B$.  Then $({\rm ker}\phi , \leq_\phi)$ is a point of
$\sper A$.  We shall simply denote this point by the letter $\phi$.
For any $a \in A$, we may identify $A[\phi]$ with the image of $A$ in
$B$ under the map $\phi$ and $a(\phi)$ with the image of  $a$.

Let $A = k[x, y]$, where  $k \subseteq \R$  is the field of real
algebraic numbers and $\R$ is the field of real numbers.  We shall
give examples of points $\phi$ and $\psi$ in $\sper A$ both with zero
support, both centered at $(x, y)$ and both with the same value group
but with $\dim \phi =0$ and $\dim \psi = 1$.  Let  $B$ denote
$\R[[t]]$, ordered so that $0<t<\rho$ for all positive $\rho \in \R$.
The point $\phi$ is defined by setting $x(\phi) = t$ and $y(\phi) =
e^t - 1$; $\psi$ is defined by setting $x(\psi) = t$ and $y(\psi) =
\pi t$.  It is easy to see that $\Gamma_{v_{\phi}} = \Z =
\Gamma_{v_{\psi}}$ and that $v_\phi(a) = \ord_t(a(\phi))$ and
$v_\psi(a) = \ord_t(a(\psi))$.  We have  $\supp\phi = \{ 0\} =
\supp\psi $.  Now $K_\phi = k$ and so $\dim \phi = 0$.  In contrast,
$K_\psi = k(\pi)$ and $\dim \psi = 1$.

\bigskip
{\bf 2. Separating ideals.}  The main result of this section is
Proposition 2.2, which states a sufficient condition for the
simplicity of a separating ideal.  Before turning to this, we recall
the definition and chief properties of separating ideals and some
facts about simple ideals and successor ideals.

Let $\a ,\b\in \sper A$.  We say $a \in A$ {\it changes sign between
$\a$ and  $\b$} if either: i) $a(\a)\geq 0$ and $a(\b)\leq 0$, or ii)
$a(\a)\leq 0$ and $a(\b)\geq 0$.  The {\it separating ideal}
determined by $\a$ and $\b$, denoted $\sepab$, is the ideal of $A$
generated by the elements of $A$ which change sign between $\a$ and
$\b$.  In  [M], the following is shown:

\smallskip
\myproclaim Proposition 2.1:  If $A$ is any ring and $\a ,\b\in \sper
A$, then:

\item{a)}$\sepab$ is both an $\a$-ideal and a $\b$-ideal.The
orderings induced on $A/\sepab$ by $\a$ and $\b$ are the same (and
hence, the set of $\a$-ideals containing $\sepab$ is equal to the set
of $\b$-ideals containing $\sepab$).

\item{b)}$\sepab$ is the smallest ideal of $A$ with the properties in
(1).
\item{c)}$\sqrt{\sepab}$ is the support of the least common
specialization of $\a$ and $\b$ in $\sper A$; $A=\sepab$ if and only
if $\a$ and $\b$ have no proper common specialization.
\qed

We say that an ideal is {\it simple} if it is proper and cannot be
expressed as a product of proper ideals.  Using the notation at the
beginning of section 1, if  $J \subseteq A$ is any ideal, then
$\overline J: = \{\, b \in A \mid \exists a \in J \enskip {\it v}(b)
\geq {\it v}(a) \,\}$ is a {\it v}-ideal.  This is proper if and only
if $J$ is contained in the center of $v$.  If $I$ is a {\it v}-ideal
and $I=J_1J_2$, then also $I=\overline{J_1} \,\overline{J_2}$.  It
follows that if $A$ is local and the center of $v$ is  the maximal
ideal of $A$, then a {\it v}-ideal $I\subseteq A$ is simple if and
only if it is proper and cannot be expressed as a product of proper
$v$-ideals.

If $I$ is an $\a$-ideal (in a noetherian ring), its successor (with
respect to the valuation $v_\a$) is denoted $I^{\a}$.  This is the
largest $\a$-ideal properly contained in $I$.  In the proof of
proposition 2.2 which we are about to give, we shall consider the
successors of  $\sepab$ with respect to the orderings $\a$ and $\b$.
It is possible for $\sepab ^\a$ to equal $\sepab ^\b$ or for the two
ideals to be distinct.  If, for example, $v_\a = v_\b$ then the
successors are equal, but equality may occur even when $v_\a$ and
$v_\b$ are distinct.  Examples are given at the end of this section.

\smallskip
\myproclaim Proposition 2.2:

Let $A$ be a noetherian ring, and suppose $\a ,\b\in \sper A$ are
centered at a maximal ideal $\m \subseteq A$ and $\sepab \subseteq
\m$.  Let $k:= A/\m$.  If
$$I/I^{\a} \cong k \hbox{ for all }\a\hbox{-ideals }I \hbox{ properly
containing }\sepab, \eqno (*)$$

then $\sepab$ is simple.

\smallskip
{\it Proof}.

Pick $x\in \sepab$ satisfying the following conditions:
\itemitem{{\it i}\thinspace )} $x(\a )\geq 0$ and $x(\b)\leq 0$,
\itemitem{{\it ii}\thinspace )} ${\it v}_\a (x) = {\it v}_\a
(\sepab)$,
\itemitem{{\it iii}\thinspace )} ${\it v}_\b (x) = {\it v}_\b
(\sepab)$.

\noindent Such an $x$ exists because surely we can find $x_\a$
satisfying ({\it i}\thinspace ) and ({\it ii}\thinspace ) and $x_\b$
satisfying ({\it i}\thinspace ) and ({\it iii}\thinspace ).  Then
$x=x_\a + x_\b$ satisfies all three conditions.
Suppose $\sepab=GH$.  Without loss of generality, $G$ and $H$ are
$\a$-ideals.  Since they contain $\sepab$, $G$ and $H$ are also
$\b$-ideals.  Now it is possible to write
$x = {\sum_i g_ih_i} + e$
where $g_i \in G \setminus G^\a = G \setminus G^\b$, $h_i \in H
\setminus H^\a = H \setminus H^\b$ and $e \in \sepab^\a \cap
\sepab^\b$.  Pick any $g \in G \setminus G^\a$ and $h \in H \setminus
H^\a$.  By ($\ast$), $g_i = u_ig + g_i'$ and $h_i = v_ih + h_i'$ for
some $u_i, v_i \notin \m$, $g_i' \in G^\a$, $h_i' \in H^\a$.  Then we
have

$$x = {\sum_i (u_ig+g_i')(v_ih+h_i')} + e = gh\sum_i u_iv_i + E$$

where $v_\a (E) > v_\a (\sepab)$ and $v_\b (E) > v_\b (\sepab)$.
(This implies that $\sum_i u_iv_i \notin \m$.)  We can conclude that
$(x-E)(\a ) > 0$ and $(x-E)(\b ) < 0$. Hence, one of the factors
$\sum_i u_iv_i$, $g$ or $h$ must change sign between $\a$ and $\b$.
But this is impossible, because none of these elements is in
$\sepab$.
\qed

Under the hypotheses in the first sentence of the proposition, $k
\subseteq K_\a$ and $k \subseteq K_\b$.  (Recall that $K_\a$ is the
residue field of the valuation associated with $\a$.)  Sometimes the
containments are proper.  If $I$ is an $\a$-ideal properly containing
$\sepab $, then $I/I^\a = I/I^\b$ is a finitely generated
sub-$k$-vector space of $ K_\a$ and of $K_\b$.  Obviously, therefore,
if $k$ is real closed and either $\td(K_\a/k) = 0$ or $\td(K_\b/k) =
0$, then ($\ast$) is satisfied.  Note that  $k$ is cofinal in $K_\a$.
{\it Hence, if $k$ is the field of real numbers then ($\ast$) is
satisfied without exception.}  Indeed, we have the

\smallskip
\myproclaim Corollary 2.3:  If $A$ is the real coordinate ring of a
compact algebraic subset of $\R^n$ (e.g. $A = \R[x, y,
z]/(x^2+y^2+z^2-1)$) then every proper separating ideal defined by
points of $\sper A$ is simple.

\smallskip
{\it Proof}.  In this case, every point of $\sper A$ is centered at a
maximal ideal.  Points centered at different maximal ideals have the
unit ideal as separating ideal.
\qed

{\bf Examples 2.}  Define $V(\a, \b):=\sepab /(\sepab ^\a \cap
\sepab^\b)$.  Under the assumptions occurring in proposition 2.2,
$V(\a, \b)$ is a vector space over $k$. The residues in this space of
the elements of $A$ satisfying $a(\a) \geq 0$ and $a(\b) \leq 0$ make
up a convex cone.  We shall give some concrete illustrations.  Let $A
= \R[x, y]$.  We consider points of $\sper A$ induced by
homomorphisms $\a,\, \b : A \to \R[[t]]$, where the latter ring is
ordered so that $0<t<\rho$ for all positive $\rho \in \R$.  The image
of $a$ under the map corresponding to $\a$ is denoted $a(\a)$, and
$v_\a(a) = \ord_t(a(\a))$.

{\bf a.}  Suppose $\a,\, \b : A \to \R[[t]]$ are induced by letting
$x(\a) = t^2$,  $y(\a) = t^4 + 2t^5$, $x(\b) = t^2$ and $y(\b) = t^4
- t^5$.  Let $f = y-x^2$.  Then $f(\a)=2t^5$ and $f(\b)=-t^5$.  It is
easy to check that no polynomial $g$ with $v_\a(g)<5$ changes sign,
so $v_\a(\sepab)=5$, and $\sepab = (f, xy, y^2)$.  Moreover, if $g$
changes sign and $v_\a(g)=5$,  then $g=f+h$, where $v_\a(h)>5$, as
one can easily check.  Hence $\sepab ^\a = (x^3, xy, y^2)$ and $V(\a,
\b)$ is one-dimensional.  In this example $\sepab ^\a = \sepab ^\b$

{\bf b.}  Now let $\g$ be induced by letting $x(\g) = t^2$ and $y(\g)
= t^4 + t^5$.  Let $f_\lambda = (y-x^2)^2-\lambda x^5$.  One may
easily check that $f_\lambda (\a) \geq 0$ and $f_\lambda(\g) \leq 0$
if and only if $\lambda \in [1, 4]$.  In this case,
$v_\a(\sepag)=10$; $\sepag ^\a$ contains $f_4$, but $\sepag ^\g$ does
not.  $V(\a, \g)$ is spanned by the residues of $f_1$ and $f_4$.  The
cone $\{\,\overline{a}\in V(\a, \g)\mid a\in \sepag,\; a(\a)\geq 0,\;
a(\g)\leq 0\,\}$ is equal to the set of non-negative-linear
combinations of $\overline{f_1}$ and $\overline{f_4}$, {\it i.e.,\/},
$\{\, \mu\overline{f_\lambda}\mid \mu, \lambda\in \R,\;0\leq\mu,\;
1\leq\lambda\leq4\,\}$.

\bigskip
{\bf 3. Quadratic transforms.} We recall some information about
quadratic transforms.  Let $A$ be a regular local (noetherian) domain
with maximal ideal $\m$ and residue field $k$.  A {\it quadratic
transform} of $A$ is a local ring $B=(A[x^{-1}\m])_{\bf p}$ where
$\ord_A(x) = 1$ (here, and below, $\ord_A(x)= {\rm max}\{\,n \mid
x\in \m^n\,\}$) and ${\bf p}$ is a prime ideal of $A[x^{-1}\m]$ such
that $\m \subseteq {\bf p}$.  If $B$ is a quadratic transform of $A$,
we write $A\prec B$.   (Geometrically, if $A$ is the local ring of a
non-singular point $M$ on a variety $X$ and $\widetilde X \to X$ is
the blow-up with center $M$, then $B$ is the local ring of a point
$P$ (possibly not closed) in the inverse image of $\{ M\}$.)  Of
course, if $A$ is one-dimensional, than $B=A$.

Let $v$ be a valuation on  $A$.  We say that $v$ {\it dominates} $A$
if $v$ is non-negative on $A$ and $\cntr_A v = \m$.  If $v$ is
non-trivial (recall that this means $\cntr_A v \not= \supp_A v$),
then  the {\it transform of  $A$  along $v$} is defined to be the
ring $B = S^{-1}A[x^{-1}\m]$, where $x$ is an element of $\m$ of
minimum value and $S = \{\,a\in A[x^{-1}\m] \mid v(a) = 0 \,\}$.  It
can be shown that $B$ is a regular local domain independent of the
choice of $x$ and that $v$ dominates $B$.  Hence the process may be
repeated indefinitely, and it continues to produce {\it proper}
extensions until a discrete valuation ring is obtained, if this ever
occurs.   We frequently use the notation $A=A^{(0)} \prec
A^{(1)}\prec\cdots$ to denote the sequence of transforms along a
given valuation.  For more details, see [A].  If $\a\in \sper A$ is
centered at $\m$ and $\supp \a \not=\m$, then the quadratic
transformation along $\a$ is defined by the valuation $v_\a$.

Transforms of ideals in $A$ are defined as follows.  Suppose $I
\subseteq A$ is an ideal with $\ord_A(I)=r$.  If $a \in I$, then
$x^{-r}a \in A[x^{-1}\m]$.  Hence, $IA[x^{-1}\m]=x^rI'$ for some
ideal $I'\subseteq A[x^{-1}\m]$.  This ideal is called the {\it
transform of $I$ in $A[x^{-1}\m]$} and $I'B$ is called the {\it
transform of $I$ in $B$}.  The transform of $I$ in $B$ is also
denoted $T_B(I)$, or $T(B)$ if $B$ is clear from context.

Certain facts regarding the real spectrum of a quadratic transform of
a local ring will be needed.  Suppose $A\subseteq B$ are any rings.
There is a functorial map $\pi :\sper B \to \sper A$, where
$\pi(\g):= \g\cap A$ and $\supp \pi (\g) = A \cap \supp \g$.  In
general, this is neither injective nor surjective, even when $B$ is a
quadratic transform of $A$.

\smallskip
\myproclaim Lemma 3.1:  Let $A$ be a regular local domain with
maximal ideal $\m$, and let $\a\in \sper A$.  Let
$B=(A[x^{-1}\m])_\p$ be a quadratic transform of $A$.  If $\supp\a
\subseteq \p\cap A$ and $\supp \a \not= \m$, then there is a unique
$\g\in \sper B$ with $\g\cap A = \a$.

\smallskip
{\it Proof.}  By standard commutative algebra, if $\q \subseteq
\p\cap A$ is any prime ideal in $A$ and $\q\not=\m$, then there is a
unique prime $\p'\subseteq B$ such that $\q = \p' \cap A$.  Note that
$x\not\in \p'$, and thus $A_\q = B_{\p'}$.  This applies in
particular to $\q = \supp \a$, provided that $\a$ satisfies the
hypotheses of the lemma.  If $A_{\supp \a} = B_{\p'}$ then the
fraction fields of the domains $A/\supp \a\subseteq B/\p'$ are equal.
Because the orderings of a domain are in one-to-one correspondence
with the orderings of its fraction field, this shows that any $\a$
satisfying the hypotheses of the lemma has a unique lift to $B$.
\qed

If $\a\in \sper A$ satisfies the conditions of the lemma, we say that
{\it $\a$ lifts to $B$.}  The unique preimage of $\a$ in $\sper B$,
when such exists, is denoted $T_B(\a)$, or simply $T(\a)$ if $B$ is
clear from context.

\smallskip
\myproclaim Lemma 3.2:  Suppose $\a,\b \in \sper A$ are non-trivial
and centered at $\m$ and that $\sepab$ is properly contained in $\m$.
Let $B$ be the first quadratic transform of $A$ along $\a$.  Then:
\item{a)}$\a$ and $\b$ both lift to $B$,
\item{b)}$T(\sepab) \subseteq \langle T(\a), T(\b) \rangle$.

\item{c)}if $A$ is of dimension 3, the containment in b) may be
proper.

\smallskip
{\it Proof.}  Because  $\sepab \not= \m$, neither $\supp\a$ nor
$\supp\b$ is $\m$.  Moreover, an element of $A$ has minimal non-zero
$v_\a$-value if and only if it has minimal non-zero $v_\b$-value.
Pick $x \in A$ with minimal non-zero value with respect to both
valuations, so $B = S^{-1}A[x^{-1}\m] $, where $S = \{\,a \in
A[x^{-1}\m] \mid v_\a(a)= 0\,\}$. By 3.1, to prove (a) it is enough
to note that $v_\a(a) \not= 0$ for all $a \in \supp \b$; this is
clear, since $\supp \b$ contains no units.

Note that $\a$ and $\b$ both lift to points of $\sper A[x^{-1}\m]$,
call them $\a'$ and $\b'$.  (Elements of $A[x^{-1}\m]$ are of the
form $f/x^n$, where $f \in \m^n$.)  To prove (b), let $\sepab '$
denote the transform of $\sepab$ in $A[x^{-1}\m]$. Then $\sepab'$ is
generated by elements of $A[x^{-1}\m]$ of the form $x^{-n}a$, where
$a(\a) \geq 0$ and $a(\b) \leq 0$.  Any such element belongs to
$\langle \a',\b' \rangle$, since $x$ does not change sign between
$\a$ and $\b$.  Hence, $\sepab' \subseteq \langle \a',\b' \rangle$,
so $T(\sepab) = \sepab' B \subseteq \langle T(\a), T(\b) \rangle$.

The example needed for the last assertion is given below (3.e).
\qed

The proof just given shows, in fact, that if $B=(A[x^{-1}\m])_{\bf
p}$ is any first quadratic transform of $A$ and if $x \not\in
\sepab$, then $T(\sepab) \subseteq \langle T(\a), T(\b) \rangle$
provided the lifts $T(\a)$ and $T(\b)$ exist.

In dimension two, the transforms of separating ideals behave very
well, as we show in the next section.

\smallskip
{\bf Examples 3.}  In these examples, we use $A'$ to denote a
quadratic transform of $A$.  Also, we write $I'$ in place of $T(I)$
to denote the ideal transform, $\a'$ to denote the lift of $\a$, and
so forth.  This saves space and is easier to read.

{\bf a.}  Let $A$ be the localization of $\R[x, y]$ at the origin.
We compute the transforms $A\prec A'\prec A''\prec\cdots$ along the
order $\a$, where $\a$ is determined as in the examples at the end of
Section 2 by letting $x(\a) = t^2$ and $y(\a) = t^4 + 2t^5$.  We
shall also examine the lifts of the points $\b$ and $\g$, which, as
in the previous example, are defined by $x(\b) = t^2$, $y(\b) = t^4 -
t^5$ and $x(\g) = t^2$, $y(\g) = t^4 + t^5$.

If we let $x':=x$ and $y': = y/x$, then $A'$ is the localization of
$\R[x', y']$ at the maximal ideal $(x', y')$.  Then $\a'$, $\b'$, and
$\g'$ are defined by

$x'(\a') = x'(\b') =x'(\g') = t^2$,  $y'(\a') = t^2+2t^3$, $y'(\b') =
t^2-t^3$ and $y'(\g') = t^2+t^3$.  It is easy to check that
$v_\a(\langle \a', \b'\rangle) = 3$ and hence $\langle \a',
\b'\rangle = (y'-x', {x'}^2) = \sepab'$.  Also, $v_\a(\langle \a',
\g'\rangle) = 6$, and hence  $\langle \a', \g'\rangle= ((y'-x')^2,
{x'}^3, {x'}^2{y'})= \sepag'$.

The transform $A''$ of $A'$ along $\a'$ is the localization of
$\R[x'', y'']$ at $(x'', y'')$  where $x'':= x'$ and $y'':= y'/x' -
1$.  The lifts of $\a'$, $\b'$ and $\g'$ are defined by  $x''(\a'') =
x''(\b'') =x''(\g'') = t^2$, and  $y''(\a'') = 2t$,

$y''(\b'') = -t$ and
$y''(\g'') = t$.  Now,

$v_\a(\langle \a'', \b''\rangle) = 1$ and hence $\langle \a'',
\b''\rangle = (x'', y'') = \sepab''$.  Also, $v_\a(\langle \a'',
\g''\rangle) = 2$ and from this one easily deduces $\langle \a'',
\g''\rangle = (y'', {x''}^2)= \sepag''$.

$A'''$ is the localization of $\R[x''', y''']$ at $(x''', y''')$
where $x''':= x''/y''$ and $y''':= y''$.  One shows easily that
$\langle \a''', \g'''\rangle= (x''', y''') = \sepag'''$.

{\bf b.}  We illustrate the remark about the non-uniqueness of points
of  $\sper A'$ which contract to a given point of $\sper A$ with
support $\m$.  Let $\delta^+,\delta^- \in \sper A'$ be given by
$x'(\delta^+) =x'(\delta^-) = 0$, $y'(\delta^+) = t$ and
$y'(\delta^-) = -t$.  Then $\pi(\delta^+) = \pi(\delta^-) =$ the
unique point of $\sper A$ supported at $\m$.

{\bf c.}  An example in which $A$ is a two-dimensional regular local
ring and $A'$ is a discrete valuation ring may be obtained as
follows.  Let $k$ be the real algebraic numbers, and let $A$ be the
localization of $k[x, y]$ at the origin.  Suppose $\a$ is determined
by a homomorphism to $\R[[t]]$ (ordered as in the example at the end
of Section 1 ) with $x(\a) = t$ and $y(\a) = \pi t$.  In this case,
$v_\a = \ord_A$; $A'$ is the localization of $k(\pi)[x]$ at $x=0$.

{\bf d.}  If $\sepab = \m$, it is possible that $\langle \a', \b'
\rangle \subseteq \m_1$ (the maximal ideal of $A'$), showing that
$\sepab' \not\subseteq \langle \a', \b' \rangle$ is possible when the
hypothesis in 3.2 is not met.  Let $A$ be the localization of $\R[x,
y]$ at the origin and let $\a$ and $\b$ be determined by letting
$x(\a) = t$, $y(\a) = t^2$, $x(\b) = -t$ and $y(\b) = -t^2$.  Then
$\sepab = \m$.  As $x'(\a') = t$, $y'(\a') = t$, $x'(\b') = -t$ and
$y'(\b') = t$, we see that $\langle \a', \b' \rangle = \m_1$.  The
points $\a$ and $\b$ may be thought of as having a common tangent,
but as determining different directions along that tangent, which
explains geometrically why these points retain a common center after
a quadratic transform.  These points do not retain a common center
after another transform, however.

{\bf e.}  The following example was found using a computer program by
Alvis which searches systematically for polynomials which change sign
between given orders.  The details will appear elsewhere.  We take
$A$ to be the localization of $\R[x, y, z]$ at the origin.  For any
real number $u$, let $\g_u$ be determined by letting

$$x(\g_u) = t^6$$
$$y(\g_u) = t^{10} + ut^{11}$$
$$z(\g_u) = t^{14} + t^{15}.$$
We examine the ideal  $I : =  \langle \g_1, \g_3 \rangle$.   Let $f
:= z^2 - x^3y$.  Then $v_{\g_1}(f) = 29$.  Exhaustive search shows
that $f$ has minimal $v_{\g_1}$-value among polynomials which change
sign between  $\g_1$ and $\g_3$.  Knowing this, it is routine to find
a list of generators for $I$; we find $I=(f, x^5, x^4y, x^3z, x^2y^2,
xyz, xz^2, y^3, y^2z, yz^2, z^3)$ (with some obvious redundancy in
the list of generators).  Clearly $\ord_A(I) = 2$.  Let  $x' = x$,
$y' = y/x$ and  $z' = z/x$.  The quadratic transform of $A$ along
$\g_1$, $A'$, is the localization of $\R[x', y', z']$ at the origin.
The transform $I'$ of $I$ is not a valuation ideal, so it is clearly
not a separating ideal.  But the situation is even more complicated.
If $v$ is the valuation induced by $\g_1$, we have $v(I')=17$, but
$v(\langle \g_1', \g_3' \rangle) = 13$, since $2y'z' - x'^2 - y'^3$
changes sign between  $\g_1'$ and  $\g_3'$.  (This shows that it is
not possible to obtain an equality by using the ``complete
transform'' defined in [L] in place of the ideal transform we have
used.  The strange behavior of this ideal seems to be related to the
fact that it is not finitely supported---see [L] for the meaning of
this term.)

\bigskip
{\bf 4. Two-dimensional rings}.  We begin this section by recalling
the results from Zariski's theory of valuation ideals in
two-dimensional regular local rings which we will be using.  We make
no attempt to give any indication of the algebra required to derive
these and refer the reader to [Z], [ZS] or [Hu] for details.

In this section, we make the following assumptions:
\item{$i)$}$A = (A, \m, k)$ is a two-dimensional regular local
domain.  (Zariski's theory requires no special assumptions about the
residue field $k$, but for our purposes we shall eventually need to
assume $k$ is real closed.)

\item{$ii)$}$v$ is a valuation centered at $\m$ which is non-trivial
({\it i.e.}, $\cntr v \not= \supp v$).

The following notational conventions are used in this section.
Assume that $A = A^{(0)} \prec A^{(1)} \prec \dots$ is the sequence
of quadratic transforms along $v$.  The maximal ideal and residue
field of $A^{(i)}$ are denoted $\m^{(i)}$ and $k^{(i)}$,
respectively.  Let $\{\, I_i \mid i=0,1,2,\dots\,\}$ be the initial
segment of the sequence of ${\it v}$-ideals of $A$ ($I_0=A$,
$I_1=\m$).  Similarly, let $\{\, J_i \mid i=0,1,2,\dots\,\}$ be the
initial segment of the sequence of ${\it v}$-ideals of $A^{(1)}$.
Also, let $\{\,{\cal I}_i = I_{n_i} \mid i= 0, 1, \dots \,\}$ be the
subsequence of $\{ I_i\}$ consisting of the simple ${\it v}$-ideals
of $A$ (with ${\cal I}_0=\m$), and let $\{\, {\cal J}_i = J_{n_i}
\mid i=0, 1, \dots\,\}$ be defined similarly.

Let $T$ be the ideal transform operation corresponding to the passage
from $A$ to $A^{(1)}$.  The {\it inverse transform}, denoted $W$, is
defined as follows.  Suppose $B$ is a first quadratic transform of
$A$ and $J\subseteq B$ is an ideal.  Since $J$ is finitely generated,
there is an integer $n$ such that $x^nJ=IB$ for some ideal
$I\subseteq A$.  The {\it inverse transform of $J$} is the ideal
$W(J) = x^{n_0}J \cap A$, where $n_0$ is the least such integer.
Observe that $T(W(J))=J$, but in general $W(T(I)) \supseteq I$ only.
(This definition for $W$ seems to be appropriate only in dimension 2.
For a generalization, see [L], proof of 2.3.)

\smallskip
\myproclaim 4.1:  (See [ZS], p. 390.)  If $v \not= \ord_A$, then:
\item{a)}The transform in $A^{(1)}$ of any $I_i$ is a member of the
sequence $\{ J_i\}$.
\item{b)}The inverse transform of any $J_i$ is a member of the
sequence $\{ I_i\}$.
\item{c)}Any $I_i$ is of the form $\m^hJ$, where $J$ is the inverse
transform of some $J_j$.
\item{d)}If $W(J_i)\subseteq W(J_j)$, then $J_i\subseteq J_j$.
\qed

If $\td(K_v/k) = 1$, then $v$ is said to be a {\it prime divisor}.
For example, if $B$ is obtained from $A$ by a finite sequence of
quadratic transforms, then $\ord_B$ is a prime divisor, whose center
in $A$ is $\m$.

\smallskip
\myproclaim 4.2: (See [ZS], p. 391--2.) The sequence $\{{\cal I}_i\}$
is finite if and only if $v$ is a prime divisor.  If there are
exactly $n$ $\m$-primary simple $v$-ideals $\{\, {\cal I}_0, \dots ,
{\cal I}_{n-1}                         \,\}$, then the following are
true:
\item{a)}$v = \ord_{A^{(n-1)}}$.
\item{b)}$A^{(n)}$ is a discrete valuation ring and $K_v = k^{(n)}$.

\item{c)}(See [ZS], p. 363.) $k^{(i)}$  is algebraic over
$k^{(i-1)}$ for $1 \leq i \leq n-1$, while $k^{(n)}$ is a simple
transcendental extension of $k^{(n-1)}$ (indeed, $k^{(n)} =
k^{(n-1)}({s \over t})$ for any generators $s$ and $t$ of
$\m^{(n-1)})$.

\qed

In general, $I_m \supseteq I_n$ does not imply that $T(I_m) \supseteq
T(I_n)$, but the following shows that $T$ is a one-to-one
order-preserving map from the initial sequence of simple ${\it
v}$-ideals of $A$ (other than $\m$) to the initial sequence of simple
${\it v}$-ideals of $A^{(1)}$.  Let $T^{(n)}$ denote the iterated
transform along $v$ (so, if $I \subseteq A$, then $T^{(n)}(I)
\subseteq A^{(n)}$).

\smallskip
\myproclaim 4.3:  (See  Theorem 4.1 and [ZS], p. 388-9; also [Hu],
Remark 3.8.)

\item{a)} If $v \not= \ord_A$, then for $i\geq 0$, $T({\cal
I}_{i+1})= {{\cal J}_i}$ and $W({{\cal J}_i})= {\cal I}_{i+1}$.  Thus
$T^{(n)}({\cal I}_n) = \m^{(n)} \subseteq A^{(n)}$.
\item{b)} Let ${\cal P}\subseteq A$ be any $\m$-primary simple
complete ideal.  Then ${\cal P}$ uniquely determines an integer $h$
and a sequence of quadratic transforms $A=B^{(0)}\prec\cdots\prec
B^{(h)}$ such that $T_B^{(h)}({\cal P})$ is the maximal ideal of
$B^{(h)}$, and ${\cal P}$ is an $\ord_{B^{(h)}}$-ideal.  If $v$ is
any valuation centered on $A$ for which ${\cal P}$ is a $v$-ideal,
then the sequence $A=B^{(0)}\prec\cdots\prec B^{(h)}$ is the initial
part of the sequence of transforms along $v$.

\qed

To our previous assumptions,  we now add the following:
\item{$ii')$}  $\a,\b \in \sper A$,  $\cntr\a = \m = \cntr\b$ and
$\sepab$ is properly contained in $\m$.  (This ensures that $v_\a$
and $v_\b$ satisfy assumption ({\it ii\/}), above.)  $A = A^{(0)}
\prec A^{(1)} \prec \dots$ will denote the sequence of quadratic
transforms along $v_\a$.

\smallskip
\myproclaim Theorem 4.4:  Suppose $k$ is real closed and $\dim \a =
1$ (so there are only finitely many simple $\m$-primary $\a$-ideals).
If $\a \not= \b$, then $\sepab$ contains the smallest $\m$-primary
simple $\a$-ideal of $A$.

\smallskip
{\it Proof}.  Suppose the smallest $\m$-primary simple $\a$-ideal of
$A$ is ${\cal I}_{n-1}$ and suppose that ${\cal I}_{n-1} $ is {\sl
not} subset of $\sepab$.  We shall show that $\a=\b$.  Since $\sepab$
is an $\a$-ideal, $\sepab$ is a proper subset of  ${\cal I}_{n-1}$.
By 2.1.a,  each of the ideals ${\cal I}_0, \dots , {\cal I}_{n-1}$ is
a $\b$-ideal.  Hence, by 4.3.b,  $A = A^{(0)} \prec A^{(1)} \prec
\dots \prec A^{(n-1)}$ is the initial part of the sequence of
quadratic transforms along $v_\b$.  By 4.2, $v_\a =
\ord_{A^{(n-1)}}$, and $K_\a = k({s \over t})$ for any for any
generators $s$, $t$ of $\m ^{(n-1)}$.  We may choose $s$ and $t$ to
be of the form $s = {s_1 \over u}$ and  $t = {t_1 \over u}$ with
$s_1,t_1 \in {\cal I}_{n-1} \setminus ({\cal I}_{n-1})^\a$, and we
can at same time,  arrange that the transform of  $A^{(n-1)}$ along
$\b$ is a localization of $A^{(n-1)}[{s \over t}] = A^{(n-1)}[{{s_1}
\over {t_1}}]$.

Let $H_\a = \{\,h \in A \mid \res_\a h <_\a \res_\a{{s_1}\over {t_1}}
\,\}$ and $H_\b = \{\,h  \in A \mid \res_\b h <_\b \res_\b
{{s_1}\over {t_1}} \,\}$.  (Here, $\res_\a : A \to k_\a$ is the
natural map, and $\res_\b$ is defined similarly; the definition of
$k_\a$ was given just after 1.1.)   These sets, we claim, are equal.
If not---and if $h_0$ belongs to the symmetric difference---then $s_1
- h_0t_1$ changes sign between $\a$ and  $\b$, and therefore
$v_\a(s_1 - h_0t_1) \geq v_\a (\sepab)$.  Now  ${{s_1}\over
{t_1}}-h_0 \not\in \m_\a$ (since $\res_\a({{s_1}\over {t_1}}) \not\in
k$),  so $v_\a({{s_1}\over {t_1}}-h_0)=0$, and thus  $v_\a(s_1 -
h_0t_1) = v_\a(t_1) = v_\a ({\cal I}_{n-1})$.   Accordingly,
$v_\a({\cal I}_{n-1}) \geq v_\a (\sepab)$  and therefore ${\cal
I}_{n-1} \subseteq \sepab$, contrary to assumption.  Therefore, $H_\a
= H_\b$.  Since $k$ is cofinal in $K_\b$, $\res_\b {{s_1}\over
{t_1}}$ is not in $k$ and therefore $\dim \b = 1$, indeed by 4.2.c,
$v_\b = \ord_{A^{(n-1)}}$.  Since $k$ is real closed, any ordering of
a simple transcendental extension $k(\xi)$ of  $k$ is completely
determined by the set $\{\,b \in k \mid b < \xi\,\}$.  From this it
follows that $K_\a$ and $K_\b$ are not only the same field but that
the orderings induced on this field by $\a$ and $\b$ are the same.
Since $\a$ and $\b$ induce the same valuation and the induced
orderings of the residue field are identical and since the value
group is the group of integers, there are exactly two possibilities
(as shown, for instance, in [Br]): $\a =\b$ or $\b$ is induced by
choosing any generator of $\m^{(n)} = \m_\a$ in $A^{(n)} = V_\a$ and
assigning it the opposite sign  from that which it has in $\a$.  In
the latter case, since $s$ generates $m_\a$, either $s_1$ or $u$ must
change sign between $\a$ and $\b$.  But this implies that ${\cal
I}_{n-1} \subseteq \sepab$, contrary to assumption.  The only
possibility is $\a =\b$.

\qed

\myproclaim Lemma 4.5: Under hypotheses $(i)$ and $(ii')$ and
assuming that $k$ is real closed,  the condition:
$$I/I^{\a} \cong k \hbox{ for all }\a\hbox{-ideals }I \hbox{ properly
containing }\sepab, \eqno (*)$$
is satisfied.  In particular, under these hypotheses,  $\sepab$ is
simple.

\smallskip
{\it Proof}.  Since $k$ has no orderable algebraic extension, it
follows from 4.2.c that $I/I^{\alpha} \cong k$ provided that $I$
contains a simple $\a$-ideal.  This is clearly the case if $\a$ is
not a prime divisor, since then there are simple $\a$-ideals of
arbitrarily large value.  The case when $\a$ is a prime divisor is
handled by 4.4.
\qed

As in  Section 3, let $T(\a)$ denote the lift of $\a$ to a point of
$\sper A^{(1)}$.  By Lemma 3.2, $T(\b)$ exists.

\smallskip
\myproclaim Lemma 4.6: Under hypotheses $(i)$ and $(ii')$ and
assuming the condition $(*)$ of 2.2 and 4.5,  $$W(\langle T(\a),T(\b)
\rangle) \subseteq \sepab.$$

\smallskip
{\it Proof.}  Let $I = W(\langle T(\a),T(\b) \rangle)$.  If $I/I^\a
\not\cong k$, then the conclusion follows from ($\ast$).  Otherwise,
argue as follows.  Since $\sepab \not= \m$, $x$ (being of minimal
non-zero $v_{\a}$-value) does not change sign between $\a$ and $\b$.
Pick $s \in \langle T(\a),T(\b) \rangle$ which changes sign between
$T(\a)$ and $T(\b)$ and is of minimal $v_{\a}$-value with this
property.  Then (referring to the notation at the end of the
paragraph preceding 4.1), $x^{n_0}s = wu$ for some $w \in I$ and some
unit $u \in A^{(1)}$ (since $I/I^\a \cong k$).  Because neither $x$
nor $u$ changes sign between $T(\a)$ and $T(\b)$, $w$ must change
sign. Hence $w \in \sepab$.  This shows that $v_{\a}(I) =
v_{\a}(s)+n_0 \geq v_{\a}(\sepab)$, which implies $I \subseteq
\sepab$.
\qed

\myproclaim Theorem 4.7:  Suppose $A = (A, \m, k)$ is a
two-dimensional regular local domain with real closed residue field
$k$.  Also suppose $\a,\b \in \sper A$,  $\cntr\a = \m = \cntr\b$ and
$\sepab$ is properly contained in $\m$.  Then:
\item{a)}$T(\sepab) = \langle T(\a), T(\b) \rangle$,

\item{b)}$W(\langle T(\a), T(\b) \rangle) = \sepab$.

\smallskip
{\it Proof}. By Lemmas 3.2, 4.5 and 4.6,

$$T(\sepab) \subseteq \langle T(\a), T(\b) \rangle \quad {\rm and}
\quad W(\langle T(\a), T(\b) \rangle) \subseteq \sepab.$$
Our assumptions guarantee that $v_\a$ is not the order valuation.
Applying $W$ and $T$ to the containments we have and using 2.2 and
4.3, we get
$$\sepab\subseteq W(\langle T(\a), T(\b) \rangle) \quad {\rm and}
\quad \langle T(\a), T(\b) \rangle \subseteq T(\sepab).\eqno /////$$

\bigskip

{\bf 5. Applications.}  Throughout this section, $(A,\m)$ is a
2-dimensional regular local ring with real closed residue field.  As
an immediate corollary of 4.7, we have

\smallskip
\myproclaim Proposition 5.1:  An $\m$-primary ideal of $A$ is a
separating ideal if and only if it is a simple complete ideal.

\smallskip
{\it Proof.}  One direction is immediate from 4.5 and 2.2.  Suppose
that ${\cal P}$ is a simple complete ideal, and
$A=B^{(0)}\prec\cdots\prec B^{(h)}$ is the sequence of transforms in
4.3.b.  The maximal ideal of $B^{(h)}$ is a separating ideal---for
example, let $u$ and $v$ be a pair of generators for this ideal and
let $\a$ and $\b$ have as their supports the ideals $(u)$ and $(v)$,
respectively.  Note that $\a$ and $\b$ are the lifts of their
restrictions to $B^{(i)}$ for any $0\leq i < h$.  Applying 4.7
repeatedly shows that ${\cal P} = \langle\a\cap A, \b\cap A\rangle$.
\qed

Two different proofs of the Pierce-Birkhoff conjecture for smooth
surfaces defined over an arbitrary real closed field can be given
based on the results of section 4.  In [MR], a proof based on the
following approximation theorem is given:

\smallskip
\myproclaim Proposition 5.2:  Suppose $\sepab$ is $\m$-primary.  Let
$X,Y\subseteq \sper A$ be closed constructible sets with $\a\in X$
and $\b\in Y$.  Then there are points $\a'\in X$ and $\b' \in Y$ with
$\dim\supp\a' = \dim\supp\b'\leq 1$ and $\langle \a',\b'\rangle =
\sepab$.
\qed

Actually, this is not difficult to prove directly when the residue
field of $A$ is the real numbers, in which case one is assured that
$K_\a = \R$.  In general, however, $K_\a$ my be a transcendental
extension of the residue field of $A$.  [MR] first proves 5.2 in the
special case that $\sepab$ is maximal.  To prove the general case,
quadratic transforms are applied until the separating ideal is
maximal, the special case is invoked, and then the solution is
transported back to the original ring via 4.7.

In [MS], a different and more abstract approach to the
Pierce-Birkhoff conjecture for smooth surfaces is given.  Here, 4.7
is used to prove:

\smallskip
\myproclaim Proposition 5.3:  Suppose that $\phi$ is an abstract
semialgebraic function on $\sper A$ with the property that for all
$\gamma \in \sper A$, there is $a\in A$ with $\phi(\g)= a(\g)$, {\it
i.e.}, $\phi$ is an ``abstract piecewise polynomial''.  Let $\a$ and
$\b$ be any points of $\sper A$.  Suppose $\phi(\a)=f(\a)$ and
$\phi(\b)=g(\b)$  for some $f,g\in A$.  Then $f-g \in \sepab$.
\qed

In the terminology of [M], 5.3 says that $A$ is a ``Pierce-Birkhoff
ring''.  In proving 5.3, the only cases that require any work are
when $\sepab$ is $\m$-primary.  As with 5.2, this is handled applying
quadratic  transforms to simplify the problem, and then transporting
information back to the original, untransformed ring by means of 4.7.
{}From 5.3, [MS] deduces that any regular  two-dimensional algebra over
a real closed field is Pierce-Birkhoff.  This result includes as a
special case the Pierce-Birkhoff conjecture for smooth surfaces.
Another result to be presented in [MS] depends upon defining
abstractly the direction of $\a$ in the Zariski tangent space at
$\cntr \a$.  It is shown, using 4.7, that $\sepab$ is the simple
ideal corresponding to the first of the transforms along $\a$ in
which the abstract directions of $\a$ and $\b$ are distinct.

\beginsection References

\item{[A]} S.~Abhyankar, On the valuations centered in a local
domain, American Journal of Mathematics, 78(1956), 321-348.

%\item{[B]} N.~Bourbaki, Commutative algebra, Hermann, Paris, 1972.

\item{[BCR]} J.Bochnak, M.Coste and M.~F.~Roy, G\'eom\'etrie
alg\'ebrique r\'eelle, Ergebnisse der Mathematik und ihher
Grenzgebiete, 3. Folge, Band 12, Springer-Verlag,
Berlin-Heidelberg-New York, 1987.

%\item{[BKW]}A.~Bigard, K.~Keimel, and S.~Wolfenstein, Groupes et
%anneaux reticules, Lecture Notes in Math.\ {\bf608}, Springer-Verlag,
%Berlin and New York, 1977.

\item{[Br]} R.~Brown, Real places and ordered fields, Rocky Mountain
J. Math. 1(1971), 633-636.

%\item{[Bru]}G.~Brumfiel, {\sl Partially Ordered Rings and
%Semi-Algebraic Geometry\/}. London Math.\ Soc.\ Lecture Note Series
%{\bf37}, Cambridge, 1979. {\sl Zbl.\/\ \bf415} (1980), 13015. {\sl
%MR\/ \bf81e}: 55029.

%\item{[DM]} C.~N.~Delzell and J.~J.~Madden, Lattice-ordered rings
%and semialgebraic geometry: I.  {\sl Proceedings of Conference on
%Real Analytic and Algebraic Geometry II, Trento, 1992\/}. De Gruyter,
%Berlin.

\item{[Hu]} C.~Huneke, Complete ideals in two-dimensional regular
local rings, Commutative algebra, Proc. Microprogram, MSRI
publication no. 15, Springer-Verlag, Berlin-Heidelberg-New York,
1989, 325--338.

\item{[L]} J.~Lipman, Complete ideals in regular local rings,
Agebraic geometry and commutative algebra in honor of Masayoshi
Nagata, 1987, 203--231.

\item{[M]}  J.~Madden, Pierce-Birkhoff rings, Arch. Math., 53(1989),
565--570.

%\item{[Mah]}  L.~Mah\'e, On the Pierce-Birkhoff Conjecture, Rocky
%Mtn.J.Math. 14(1984), 983--985.

\item{[MS]} J.~Madden and R.~Robson, In preparation.

\item{[MS]} J.~Madden and N.~Schwartz, In preparation.

\item{[Z]} O.~Zariski, Polynomial ideals defined by infinitely near
base points, Amer. J. Math., 60(1938), 151--204.

%\item{[Zar2]} O.Zariski, The resolution of singularities of an
%algebraic surface, Ann. of Math., 40(1939), 639--689.

%\item{[Zar3]}O.~Zariski, {\sl Algebraic Surfaces\/}, 2nd ed.
%Springer-Verlag, 1971.

\item{[ZS]} O.~Zariski and P.~Samuel, Commutative algebra, vol.II,
Van Nostrand 1960.

\beginsection Authors:

Dean Alvis

Department of Mathematics and Computer Science

Indiana University at South Bend

South Bend, IN 46615-1408  USA

\smallskip

\noindent Bernard Johnston

Department of Mathematics

Florida Atlantic University

Boca Raton, FL 33431-1505  USA

\smallskip
\noindent James J. Madden

Department of Mathematics

Louisiana State University

Baton Rouge, LA 70808-6710  USA

\beginsection Suggested shortened version of title:

Complete ideals defined by sign conditions

\bye